\documentclass[11pt]{article}
\usepackage{mathptmx} %%% for ps fonts
\usepackage{amsmath}\usepackage{epsf,amsfonts,amssymb}
\usepackage{theorem}  % For \theorembodyfont to define 
\makeatletter
\@addtoreset{equation}{section}

\makeatother

\makeatletter
\def\operatorname#1{\mathop{\operator@font #1}\nolimits}%
\makeatother
\renewcommand{\P}{\Lambda}
\newcommand{\Op}{{\mbox{Op}}}
\renewcommand{\exp}{\operatorname{exp}}
\newcommand{\ad}{\operatorname{ad}}
\newcommand{\cyclic}{\mathop{\kern0.9ex{{+}
 \kern-2.15ex\raise-.25ex\hbox{\Large\hbox{$\circlearrowright$}}}}\limits}

\newcommand{\w}{\wedge}
\newcommand{\p}{\partial}
\newcommand{\raa}{\rightarrow}

\newcommand{\R}{\mathbb{R}}

\newcommand{\op}[1]{\!\!\mathop{\rm ~#1}\nolimits}

\newtheorem{theorem}{Theorem}
\newtheorem{prop}{Proposition}
\newtheorem{lemma}{Lemma}
\newtheorem{cor}{Corollary}

{\theorembodyfont{\normalfont\rmfamily}
\newtheorem{defi}{Definition}
\newtheorem{rem}{Remark}
}
\makeatletter

\newlength{\blackboxsize} \blackboxsize=1.2ex
\def\blackbox{\hbox{\vrule height  \blackboxsize 
width  \blackboxsize depth 0ex}}
\newenvironment{proof}[1][]%
 {\def\proof@temp{#1}\par\noindent
  \textsc{Proof}\ifx\proof@temp\@empty\else\ (#1)\fi\hspace{1em}}
 {\par~~\hfill\blackbox\par\vspace{.5\baselineskip}}
\makeatother

\begin{document}
%%
%% Top matter
%%
\title{Universal star products}
\author{
Mourad Ammar$^{1, 2}$ \\
\textit{\small  $1$ Universit\'e du Luxembourg} \\
\textit{\small  avenue de la fa\"\i encerie} \\
\texttt{\small mourad.ammar@uni.lu} \\
  \and\\
V\'eronique Chloup$^{ 2}$\\
\textit{\small  $2$ Universit\'e Paul Verlaine - Metz, LMAM}\\
 \textit{\small Ile du Saulcy, F-57045~Metz Cedex 01,France}\\
 \texttt{\small chloup@poncelet.univ-metz.fr}\\
    \and\\
Simone Gutt$^{3,2}$\\
\textit{\small Membre de l'Acad\'emie Royale de Belgique} \\ 
\textit{\small $3$ Universit\'e Libre de Bruxelles, Campus Plaine CP 218,} \\ 
\textit{\small Bvd du Triomphe, B-1050~Brussels, Belgium}\\
\texttt{\small sgutt@ulb.ac.be, gutt@poncelet.univ-metz.fr}
}

\date{}

\maketitle

\setcounter{page}{0}

\thispagestyle{empty}

\begin{abstract}
One defines the notion of universal deformation quantization:
given any manifold $M$, any Poisson structure $\P$ on $M$
and any torsionfree linear connection $\nabla$ on $M$,
a universal deformation quantization associates to this data a star product
on $(M,\P)$  given by a series of bidifferential operators
whose corresponding tensors are given by universal polynomial expressions
 in the Poisson tensor $\P$, the curvature
tensor $R$ and their covariant iterated derivatives.
Such universal deformation quantization exist.
We study their  unicity  at order 3
in the deformation parameter, computing
the appropriate universal Poisson cohomology.
\end{abstract}



\newpage
\section{Introduction}\label{sect:intro}
Let $(M,\P)$ be a d-dimensional Poisson manifold. 
Let $C^{\infty}(M)$ be the commutative algebra of 
$\mathbb{K}$-valued smooth function on $M$, 
where $\mathbb{K}$ is $\mathbb{R}$ or $\mathbb{C}$. 
The Poisson bracket  of functions on $(M,\P)$ is denoted by $\{f,g\}:=\P(df,dg),~ f,g\in C^{\infty}(M)$, 
and $m(f,g)=f.g$ denotes the pointwise multiplication. 
\begin{defi}\cite{BFFSL}
Let $C^{\infty}(M)[[\nu]]$ be the space of formal power series in the formal parameter 
$\nu$ with coefficients in $C^{\infty}(M)$. 
A {\bf differential star product} on $(M,\P)$ is a bilinear map  
\[
* : C^{\infty}(M)\times C^{\infty}(M)\raa C^{\infty}(M)[[\nu]]
\quad\quad(f,g)\raa f*g:=f.g+\sum_{r\geq1}\nu^rC_r(f,g),
\]
such that
\begin{enumerate}
  \item its $\mathbb{K}[[\nu]]$-bilinear extension is an associative product $ (f*g)*h=f*(g*h); $
  \item $C_1(f,g)-C_1(g,f)=2\{f,g\};$
  \item each $C_r$ is a bidifferential operator vanishing on constants.
\end{enumerate}
\end{defi}
Kontsevich gave in \cite{K} an explicit formula for a star product on $\R^d$
endowed with any Poisson structure, 
as a special case of his formality theorem for the Hochschild complex of
multidifferential  operators.
\begin{theorem}\cite{K}
The formality theorem on $\R^d$ is proven  through
the explicit formula of an $L_\infty$ morphism from the differential graded Lie algebra
of polyvectorfields on $\R^d$ to the differential graded Lie algebra of polydifferential
operators on $\R^d$.
 This
consists in  a collection of  
multilinear graded symmetric maps $U_j$ associating to $j$ multivectorfields
$\alpha_k \in \Gamma(\R^d,\wedge^{m_k}T\R^d)$
a multidifferential operator 
$$
U_j(\alpha_1,\ldots,\alpha_j ):C^\infty(\R^d)^{\otimes r} \rightarrow C^\infty(\R^d)
$$
with $r=\sum_k m_k-2j+2$, satisfying  quadratic relations $(*)$ which translate the fact
that they are the Taylor coefficients of a $L_\infty$ morphism.
The maps $U_j$ are defined in terms of graphs; in particular the coefficients
of the multidifferential operator $U_j(\alpha_1,\ldots,\alpha_j )$
are given by multilinear  universal expressions in the partial derivatives
of the coefficients of the multivectorfields $\alpha_1,\ldots,\alpha_j .$\\
The relations $(*)$ imply that, given any Poisson structure $\P$ on $\R^d$, 
the formula 
\begin{equation}
f*^\P_Kg=f\,g+\sum_{n=1}^{\infty}\frac{\nu^n}{n!}U_n(\P,\ldots,\P)(f,\,g)=fg+\nu \P(df,dg)+O(\nu^2)
\end{equation}
defines a star product on $(\R^d,\P)$. 
Observe that  each
$U_n(\P,\ldots,\P)$ is a bidifferential operator of order maximum $n$ in each argument
whose coefficients 
are   polynomials of degree $n$ in the partial derivatives of the coefficients 
of the tensor $\P$.
\end{theorem}
Kontsevich  then obtained the existence of star products on a general Poisson manifold
using abstract arguments. 

A more direct 
 construction of a star product on a $d$-dimensional
 Poisson manifold $(M,P)$, using Kontsevich's 
formality on $\R^d$, was given by Cattaneo, Felder
 and Tomassini in \cite{CFT}.
 Given a torsionfree connexion $\nabla$ on $(M,P)$
one builds an  identification of the commutative algebra 
 $C^\infty(M)$ of smooth functions on  $M$
 with the algebra  of  flat sections  of the jet bundle 
 $E\rightarrow M$, for the Grothendieck connection $D^G$.
 The next point is to  ``quantize" this situation: 
 a deformed algebra stucture on $\Gamma(M,E)[[\nu]]$
 is obtained through fiberwize quantization of the jet bundle using Kontsevich 
 star product on $\R^d$, and a deformed flat connection $D$ which is a derivation of 
 this deformed algebra structure is constructed  ``\`a la Fedosov".
Then one constructs an identification between the formal series of  functions on $M$
 and the algebra of flat sections of this quantized bundle of algebras; this identification
 defines the star product on $M$.

 Later, Dolgushev \cite{D}  gave in a similar spirit a construction for a 
 Kontsevich's formality quasi-isomorphism for a general smooth
 manifold. The construction starts again with a torsionfree linear connexion
 $\nabla$ on $M$ and the identification of the commutative algebra 
 $C^\infty(M)$ of smooth functions on  $M$
 with the algebra  of  flat sections  of the jet bundle 
 $E\rightarrow M$, for a connection $D^F$ constructed ``\`a la Fedosov".
 This is extended to a resolution of the space  ${T}_{poly}(M)$
 of polyvectors on the manifold
 using the complexes of forms on $M$ with values
 in the bundle of formal fiberwize polyvectorfields on $E$ and
a resolution of the space  ${D}_{poly}(M)$
of polydifferential operators on the manifold
 using the complexes of forms on $M$ with values
 in the bundle of formal fiberwize polydifferential operators on $E$.
  The fiberwize Kontsevich $L_\infty$ morphism   is then twisted and contracted
 to yield a $L_\infty$- morphism from  ${T}_{poly}(M)$ to 
 ${D}_{poly}(M)$.\\

Given a torsionfree linear connection $\nabla$ on a manifold $M$, 
any multidifferential operator $\Op:(C^{\infty}(M))^k\rightarrow C^{\infty}(M)$
writes in a unique way as
\begin{equation}
\Op(f_1,\ldots,f_k)=\sum_{J_1,\ldots,J_k}\Op^{J_1,\ldots,J_k}\nabla^{sym}_{J_1}f_1\ldots \nabla^{sym}_{J_k}f_k
\end{equation}
where the $J_1,\ldots,J_k$ are multiindices and  $\nabla^{sym}_{J}f$
is the symmetrised covariant derivative of order $\vert J \vert$ of $f$:
\[
\nabla^{sym}_{J}f=\sum_{\sigma\in S_m} \frac{1}{m!}\nabla^m_{i_{\sigma(1)}\ldots i_{\sigma(m)}}f
\qquad \rm{for }~J=(i_1,\ldots,i_m),
\]
where $ \nabla^m_{i_1\ldots i_m}f:=\nabla^m f(\partial_{i_1},\ldots,\partial_{i_m})$ 
with $\nabla^m f$ defined inductively by $\nabla f:=df$ and\\
$\nabla^mf(X_1,\ldots,X_m)=(\nabla_{X_1}(\nabla^{m-1}f))(X_{2},\ldots,X_{m})$.\\
The tensors $\Op^{J_1,\ldots,J_k}$ are covariant tensors of order
$\vert J_1\vert +\ldots +\vert J_k\vert$ which are symmetric within each block
of $J_r$ indices; they are called {\bf the tensors associated} to  $\Op$ for the given connection.
\begin{defi}
A {\bf universal Poisson-related multidifferential operator} will be the association to
any manifold $M$, any torsionfree connection $\nabla$ on $M$
and any Poisson tensor $P$ on $M$, of a multidifferential operator
$\Op^{(M,\nabla,\P)}:(C^{\infty}(M))^k\rightarrow C^{\infty}(M)$, so that,
the tensors associated to  $\Op^{(M,\nabla,\P)}$ for $\nabla$  are given by universal polynomials  
in $\P$,  the curvature
tensor $R$ and their covariant multiderivatives, involving concatenations.\\
We shall say that a universal Poisson-related multidifferential operator
is {\bf of no-loop type} if the concatenations only arise
between different terms, not within a given term
(i.e.  $(\nabla_r \P)^{is}(\nabla_s \P)^{jr}{\nabla}^{sym}_{ij}$ is 
of no-loop type
but $(\nabla_t \P)^{ti}{\nabla}^{sym}_{i}$ or $R^r_{str}\P^{si}\P^{tj}{\nabla}^{sym}_{ij}$
are not).\\
We shall say that the universal Poisson-related operator $\Op^{(M,\nabla,\P)}$
is  {\bf a polynomial of degree $r$ in the Poisson structure} 
if $\Op^{(M,\nabla,t\P)}=t^r\Op^{(M,\nabla,\P)}$.
\end{defi}
 \begin{defi}
A {\bf universal star product} $*=m+\sum_{r\geq 1}\nu^rC_r$ will be the association
to any manifold $M$, any torsionfree connection $\nabla$ on $M$
and any Poisson tensor $P$ on $M$, of a
differential star product $*^{(M,\nabla,\P)}:=m+\sum_{r\geq 1}\nu^rC_r^{(M,\nabla,\P)}$ where
each $C_r$ is a universal Poisson-related bidifferential operator
of no-loop type, which is a polynomial of degree $r$ in the Poisson structure.
\end{defi}
An example of a universal star product at order 3 is given 
in  section \ref{section:exampleorder3}.\\
Unicity at order $3$ is studied in section \ref{section:unicity}
using  universal Poisson cohomology. This, we   compute for 
universal Poisson-related bidifferential operators of order $1$ 
in each argument defined by low order polynomials in the Poisson structure.\\

\noindent The existence of a universal star product is implied either by the
globalisation proof of Cattaneo, Felder and Tomassini \cite{CFT,CF}
using the exponential map of a torsionfree linear connection, either by the 
globalisation of the formality given by Dolgushev \cite{D}.
We show this existence in section \ref{section:existence}, stressing
first the relations between the resolutions involved
in the two constructions in section \ref{section:CFT-D}.

\section{An example at order 3}\label{section:exampleorder3}
\begin{theorem} 
There exists  a  universal star 
product up to order three, which associates to a Poisson manifold 
$(M,\P)$ and a torsionfree linear connection $\nabla$ on $M$,
the star product at order three defined by 
\begin{equation}\label{eq:staratorder3}
f{\tilde{*}}^{(M,\nabla,\P)}_3g=f.g+\nu\{f,g\}+\nu^2{\tilde{C}}^{(M,\nabla,\P)}_2(f,g)
+\nu^3{\tilde{C}}^{(M,\nabla,\P)}_3(f,g),\quad\quad f,g\in C^{\infty}(M)
\end{equation}
for 
\begin{equation}
{\tilde{C}}^{(M,\nabla,\P)}_{2}(f,g) =\frac{1}{2}\P^{kr}\P^{ls}\nabla _{kl}^{2}f\nabla _{rs}^{2}g 
+\frac{1}{3}\P^{kr}\nabla _{r}\P^{ls}(\nabla
_{kl}^{2}f\nabla _{s}g+\nabla _{s}f\nabla _{kl}^{2}g)
+\frac{1}{6}\nabla_l\P^{kr}\nabla_k\P^{ls}\nabla_rf\nabla_sg,
\end{equation}
and
\begin{equation}
{\tilde{C}}^{(M,\nabla,\P)}_3(f,g)=\frac{1}{6}S^{(M,\P)~3}_{\nabla}(f,g)=-\frac{1}{6}\P^{ls}(
{\mathcal{L}}_{X_{f}}\nabla )_{kl}^{j}({\mathcal{L}}_{X_{g}}\nabla )_{js}^{k} \qquad{\rm with }~X_f=i(df)\P,
\end{equation}
 where ${\mathcal{L}}_{X_{f}}\nabla$ is the tensor defined by the Lie derivative of the connection 
 $\nabla$ in the direction of the Hamiltonian vector field $X_f$
\[
({\mathcal{L}}_{X_{f}}\nabla)_{kl}^j=\P^{ij}\nabla_{kli}^3f+\nabla_k\P^{ij}
\nabla_{li}^2f+\nabla_l\P^{ij}
\nabla_{ki}^2f+\nabla_{kl}^2\P^{ij}\nabla_if+R_{ikl}^j\P^{si}\nabla_sf.
\]
\end{theorem}
This can be seen by direct computation.
\begin{rem}
The  operator $S^{(M,\P)~3}_{\nabla}$ was introduced by Flato, Lichnerowicz
and Sternheimer  \cite{FLS}; 
it is a Chevalley-cocycle on $(M,P)$, i.e.
\begin{equation*}
\cyclic_{u,v,w}\left\{ S^{(M,\P)~3}_{\nabla}(u,v),w\right\} +S^{(M,\P)~3}_{\nabla}(\left\{
u,v\right\} ,w)=0,
\end{equation*}
where $\cyclic_{u,v,w}$ denotes the sum over  cyclic permutations of $u,v,w$.\\
\end{rem}
For this universal star product at order $3$, there exists
a universal  Poisson-related-differential-operator-valued $1$-form 
$D$ defined as follows:

\begin{prop} 
Given any Poisson manifold 
$(M,\P)$, any torsionfree linear connection $\nabla$ on $M$,
and any   vector field $X$ on $M$, 
the differential operator $D^{(M,\nabla,\P)}_X$ defined by 
\[
D^{(M,\nabla,\P)}_Xg=Xg-\nu^2\frac{1}{6}\P^{ls}(
{\mathcal{L}}_{X}\nabla )_{kl}^{j}({\mathcal{L}}_{X_{g}}\nabla )_{js}^{k},\quad g\in C^{\infty}(M)
\]
verifies at order $3$ in $\nu$
\[
D^{(M,\nabla,\P)}_X(f{\tilde{*}}^{(M,\nabla,\P)}_3g)
-(D^{(M,\nabla,\P)}_Xf){\tilde{*}}^{(M,\nabla,\P)}_3g
-f{\tilde{*}}^{(M,\nabla,\P)}_3(D^{(M,\nabla,\P)}_Xg) =
\frac{d}{dt}|_{t=0}f{\tilde{*}}^{(M,\nabla,\phi_{t*}^{X}\P)}_3g + O(\nu^4).
\]
where $\phi^X_{t}$ denotes the  flow of the vectorfield $X$.\\
If $X$ is a Hamiltonian vector field corresponding 
to a function $f\in C^{\infty}(M)$, then $D^{(M,\nabla,\P)}_{X_{f}}$ 
coincides  with the inner derivation at order $3$ of ${\tilde{*}}_3$ 
defined by the function $f$, i.e.
\[
D^{(M,\nabla,\P)}_{X_{f}} g =\frac{1}{2\nu}(f{\tilde{*}}^{(M,\nabla,\P)}_3 g
-g{\tilde{*}}^{(M,\nabla,\P)}_3f).
\]
\end{prop}

\section{Equivalence of universal star products -- Universal Poisson cohomology}\label{section:unicity}

\begin{lemma}
\begin{itemize}
\item Any universal star product $*=m+\sum_{r\geq 1}\nu^rC_r$ is a natural star product, i.e. each bidifferential operator $C_r$ is of order at most $r$
in each argument. 
Indeed $C_r$ is a universal Poisson-related bidifferential operator
defined by a polynomial of degree $r$ in the Poisson structure;  this implies,
in view of the Bianchi's identities for the curvature tensor,
that $C_r$ is  of order at most $r$ in each argument.
\item The universal Poisson-related bidifferential operator $C_1$ of any universal star product 
is necessarily the Poisson bracket $C_1^{(M,\nabla,\P)}=\P^{ij}\nabla_i\w\nabla_j$.
\item The Gerstenhaber bracket  $[\, ,\, ]_G$ of two universal Poisson-related multidifferential operator
of degree $k$ and $l$ in $\P$, is a universal Poisson-related multidifferential operator  of 
degree $k+l$ in $\P$.
\item If a  universal Poisson-related $p$-differential operator $C$ is a Hochschild
$p$-cocycle ( where $\p:=\op{ad}\,m=[m,.]_G$ denotes the Hochschild differential) 
then $C=A+\p B$ where $A$ a  universal Poisson-related $p$-differential operator
which is of order $1$ in each argument and is the totally skewsymmetric part of $C$, and 
where $B$ is a  universal Poisson-related $(p-1)$-differential operator.
\end{itemize}
\end{lemma}
The last point comes from the explicit formulas \cite{CG,GR1} for the tensors associated to $B$ 
in terms of those associated to $C$ when one is given a connection.
\begin{defi}
A {\bf universal Poisson $p$-cocycle} is a universal Poisson-related p-differential 
skewsymmetric operator $C$ of order $1$ in each argument which is a cocycle
for the Chevalley cohomology for the adjoint representation of $(C^{\infty}(M),\{\, ,\,\})$,
i.e. with the coboundary defined by
\begin{eqnarray*}
\delta_P C(u_1,\ldots,u_{m+1})&=&\sum_{i=1}^{m+1}(-1)^i\{\, u_i, C(u_1,\ldots {\hat{u_i}}\ldots,u_{m+1})\}\\
&~&\qquad+\sum_{i<j}(-1)^{i+j} C(\{u_i,u_j\},u_1\ldots {\hat{u_i}}\ldots {\hat{u_j}}\ldots,u_{m+1}).
\end{eqnarray*}
which can be written as a multiple of 
\[
skew[\P,C^{(M,\nabla,\P)}]_G
\]
where $skew$ indicates the skewsymmetrisation in all its arguments of an operator.\\
Equivalently,   a universal Poisson $p$-cocycle $C$ 
is defined by a universal Poisson related skewsymmetric $p$-tensor $c$
(with $C(u_1,\ldots,u_p)=c(du_1,\ldots,du_p)$) so that
\[
[\P,c^{(M,\nabla,\P)}]_{SN}=0
\]
where $[\cdot,\cdot]_{SN}$ denotes the Schouten-Nijenhuis bracket of 
skewsymmetric tensors (which is the extension, as a graded 
derivation for the exterior product $\wedge$ of the usual bracket of
vectorfields).\\
A universal Poisson $p$-cocycle $C$ 
is a {\bf universal Poisson coboundary} if there exists
a universal Poisson-related skewsymmetric $(p-1)$-differential operator 
$C$ of order $1$ in each argument so that
\[
C^{(M,\nabla,\P)}=\delta_P B^{(M,\nabla,\P)} (=skew[\P,B^{(M,\nabla,\P)}]_G);
\]
(equivalently, if there exists a universal Poisson related tensor $b$ so that
$c^{(M,\nabla,\P)}=[\P,b^{(M,\nabla,\P)}]_{SN}$).\\
The {\bf universal Poisson cohomology} $H^p$ is the quotient of the space of universal
Poisson $p$-cocycles by the space of universal Poisson $p$-coboundaries.\\
We can restrict ourselves to the space of universal Poisson 
$p$-cocycles defined by polynomials
of degree $k$ in the Poisson structures and  make the quotient by  the space
of universal Poisson coboundaries defined by
polynomials of degree $k-1$. We speak then of 
the universal Poisson $p$-cohomology of degree $k$
in the Poisson structure and we denote it by $H^p_{pol k}$.
We can further restrict ourselves to universal Poisson related tensors
(or operators of order  $1$ in each argument) of no-loop type.
\end{defi}

\begin{defi}
 If $*=m+\sum_{r\geq 1}\nu^rC_r$ is a universal star product
 and if $E=\sum_{r=2}^{\infty}\nu E_r$ is a formal series  of universal 
 differential operators vanishing on constants, of no-loop type, with
each $E_r$ a polynomial of degree $r$ in the Poisson structure,
then the series $*'$ defined by 
\[
*'=(\exp\ad E )*
\]
where $\ad E \cdot=[E,\cdot]_G$,
(i.e. $f*'g=\exp E((\exp -E)f*(\exp -E)g)$),
is an equivalent universal star product.
We say that $*$ and $*'$ are universally equivalent.
\end{defi}

\begin{lemma}
If $*$ and $*'$ are universal star products
which coincide at order $k$ in the deformation parameter $\nu$,
then, by the associativity relation at order $k$,
$C'_k-C_k$ is a universal Hochschild $2$-cocycle
of no-loop type which is a polynomial of degree $k$ in the Poisson structure.
Furthermore, associativity at order $k+1$ implies
that its skewsymmetric part $p_2$ is a universal
Poisson $2$-cocycle :
\begin{equation*}
\cyclic_{u,v,w}\left\{ p^{M,\nabla,\P}_{2}(u,v),w\right\} +p^{M,\nabla,\P}_{2}(\left\{
u,v\right\} ,w)=0,
\end{equation*}
where $\cyclic_{u,v,w}$ denotes the sum over  cyclic permutations of $u,v,w$.

If it is a universal Poisson $2$-coboundary of no-loop type, then
there is a a formal series  $E$ of universal 
 differential operators vanishing on constants such that
 $(\exp\ad E )*$ and $*'$ coincide at order $k+1$.
 
 In particular, two universal star products are universally equivalent
 if $H^2_{(no-loop, pol)}=\{ 0\}$. They are always equivalent at order
 $k$ in the deformation parameter $\nu$ if $H^2_{(no-loop) pol  j}=\{ 0\}
  ~\forall 1\le j\le k.$
\end{lemma}
Consider now any  universal star product $*=m+\sum_{r\geq 1}\nu^rC_r$.
We automatically have that $C_1$ is the Poisson bracket.
Associativity at order $2$ yields $\p C_2 =\p {\tilde{C}}^{(M,\nabla,\P)}_2$ 
so
\[
C_{2}(f,g) ={\tilde{C}}^{(M,\nabla,\P)}_2 (f,g)+p_2(f,g) +\p E_2 (f,g)
\]
and the skewsymmetric part of associativity at order $3$ yields that $p_2$
is a  universal Poisson $2$-cocycle (which  is  a  
polynomial of degree $2$ in the Poisson structure).
\begin{prop}
The spaces $H^2_{pol 2}(\P)$ and $H^2_{(no-loop)pol 2}(\P)$ of universal Poisson 
$2$-cohomology of degree $2$
in the Poisson structure  vanish.
\end{prop}
\begin{proof}
The universal skewsymmetric $2$-tensors of degree 2 in $\P$ are combinations of
\begin{eqnarray*}
  && \nabla_{s}\P^{ir}\nabla_{r}\P^{js}\nabla_{i}\w\nabla_{j},\\
  &&(\nabla^2_{rs}\P^{ir}\P^{js}-\nabla^2_{rs}\P^{jr}\P^{is})\nabla_{i}\w\nabla_{j},\\
  &&  (\P^{ir}\P^{st}R_{rst}^{j}-\P^{jr}\P^{st}R_{rst}^{i})
        \nabla_{i}\w\nabla_{j},\\
  &&\P^{ri}\P^{sj}R_{rst}^{t}\nabla _i\w\nabla_j.\\
\end{eqnarray*}
The only universal cocycles are the multiples of 
\[
\P^{ri}\P^{sj}R_{rst}^{t}\nabla _i\w\nabla_j
\]
and those are the boundaries of the multiples of $\nabla_{r}\P^{ir}\p_i$.
Remark that there are no cocycles of no-loop type.
\end{proof}
Thus, universal star product at order $2$ are
 unique modulo equivalence and one can assume that $C_2={\tilde{C}}^{(M,\nabla,\P)}_2$.  Then  
the skewsymmetric part of the Hochschild $2$-cocycle $C_3-{\tilde{C}}^{(M,\nabla,\P)}_3$ is
a universal Poisson $2$-cocycle of no-loop type which is 
a polynomial of degree $3$ in the Poisson structure.
\begin{prop}
The space $H^2_{(no-loop)pol 3}(\P)$ of universal Poisson $2$-cohomology
of no-loop type and of degree $3$
in the Poisson structure  vanishes.
\end{prop}

\begin{proof}
We consider all possible universal Poisson $2$-cochains of no loop type which are polynomials of degree $3$ in the Poisson structure. They are defined by universal skewsymmetric $2$-tensors of degree 3 in $\P$ which are combinations  with constant coefficients
of the different concatenations (with no loops) of
\begin{eqnarray*}
&&\P^{\cdot \cdot}\P^{\cdot \cdot}\P^{\cdot \cdot} (\nabla^2_{\cdot \cdot}R)^{\cdot}_{\cdot \cdot \cdot} \qquad
\P^{\cdot \cdot}\P^{\cdot \cdot}\P^{\cdot \cdot}R^{\cdot}_{\cdot \cdot \cdot} R^{\cdot}_{\cdot \cdot \cdot} ,\cr
&&\P^{\cdot \cdot} (\nabla_{ \cdot}\P)^{\cdot \cdot} (\nabla_{ \cdot}\P)^{\cdot \cdot}R^{\cdot}_{\cdot \cdot \cdot}\qquad
\P^{\cdot \cdot} \P^{\cdot \cdot} (\nabla_{ \cdot}\P)^{\cdot \cdot}   (\nabla_{ \cdot} R)^{\cdot } _{\cdot \cdot \cdot}\cr
&&\P^{\cdot \cdot}\P^{\cdot \cdot} (\nabla^2_{\cdot \cdot}\P)^{\cdot \cdot} R^{\cdot}_{\cdot \cdot \cdot} \qquad
\P^{\cdot \cdot} (\nabla^2_{\cdot \cdot}\P)^{\cdot \cdot} (\nabla^2_{\cdot \cdot}\P)^{\cdot \cdot}\qquad
(\nabla^2_{\cdot \cdot}\P)^{\cdot \cdot}(\nabla_{ \cdot}\P)^{\cdot \cdot}(\nabla_{ \cdot}\P)^{\cdot \cdot}.
\end{eqnarray*}
Using the symmetry properties of $R$, the Bianchi's identities and the fact that 
$\P$ is a Poisson tensor, one is left with a combination with constant coefficients of $49$
independant terms. \\
Universal $2$-coboundaries come from the boundaries of universal $1$-tensors of degree $2$ in $\P$;
such $1$-tensors are given by combinations with constant coefficients of concatenations of
$$
 \P^{\cdot \cdot}\P^{\cdot \cdot}(\nabla_{ \cdot} R)^{\cdot } _{\cdot \cdot \cdot} \qquad
 \P^{\cdot \cdot}  (\nabla_{ \cdot}\P)^{\cdot \cdot}  R^{\cdot}_{\cdot \cdot \cdot} .
$$
Hence, modulo universal coboundaries, one can assume that the coefficients
of $4$ of the $49$ terms in a universal cochain are zero.\\
The cohomology that we are looking for is then given by the combinations with constant coefficients
of the remaining $45$ terms which are cocycles (for all possible choices of manifold, 
Poisson structure $\P$
and  connection $\nabla$.)\\
Let $C$ be a combination of those $45$ terms.  The cocycle condition is $[\P,C]_{SN}=0$.
We plug in  examples of Poisson structures and connections 
and impose this cocycle condition. This shows that all $45$
coefficients must vanish. \\
It is enough, for instance, to 
consider the example on $\R^4$, with the non vanishing coefficients of the connection defined by
\begin{eqnarray*}
&&\Gamma_{12}^1=x_1^3, \quad \Gamma_{14}^1=x_4
\quad\Gamma_{11}^2=x_1^2,\quad\Gamma_{13}^2=1,\quad\Gamma_{22}^2=1,
\quad\Gamma_{14}^2=x_3,\quad
\Gamma_{13}^3=-x_4,\\
&&\Gamma_{33}^3=1,\quad\Gamma_{44}^3=-x_2x_3x_4,
\quad\Gamma_{11}^4=1,\quad\Gamma_{13}^4=1,
\quad\Gamma_{22}^4=x_1,\quad\Gamma_{44}^3=-3.
\end{eqnarray*}
and the quadratic Poisson structure  defined by
\begin{equation*}
\P=\sum_{1=i<j=4}x_ix_j\frac{\p}{\p x_i}\w\frac{\p}{\p x_j}
\end{equation*}
From this example, one gets that $41$ of the $45$ coefficients have to vanish.
One is left with a combination with constant coefficients of four terms
and an example with constant Poisson structure in dimension $7$ shows that
all those coefficients must vanish. The non vanishing coefficients of this example are:
\begin{eqnarray*}
&&\P^{12}=1\quad\P^{15}=1\quad\P^{17}=1\quad\P^{25}=1\quad\P^{26}=2 \quad\P^{27}=2\\ &&  \P^{34}=1\quad \P^{37}=3\quad \P^{46}=1\quad\P^{47}=4\quad\P^{56}=1\quad\P^{57}=5\quad\P^{67}=6.
\end{eqnarray*}
\begin{eqnarray*}
&&\Gamma_{12}^1=1\quad\Gamma_{62}^1=x_7\quad\Gamma_{77}^1=-1 \quad\Gamma_{17}^2=x_1\quad\Gamma_{13}^3=x_6\quad\Gamma_{11}^4=1\quad\Gamma_{22}^4=1\quad\Gamma_{33}^4=2\\ && \Gamma_{77}^4=3\quad\Gamma_{12}^5=x_1x_5\quad\Gamma_{33}^5=x_2\quad\Gamma_{11}^6=1 \quad\Gamma_{44}^6=x_2 \quad\Gamma_{44}^6=x_5\quad\Gamma_{11}^7=x_7 \quad\Gamma_{44}^7=x_3\quad\Gamma_{17}^7=x_1.
\end{eqnarray*}
\end{proof}
\begin{cor}
Any universal star product is universally equivalent to one whose expression at order 
$3$ is given by formula (\ref{eq:staratorder3}).
\end{cor}

\section{Grothendieck- and Dolgushev-resolution of the space of functions}\label{section:CFT-D}

Our purpose in this section is to prove that the Fedosov-resolution 
of the algebra of smooth functions constructed in Dolgushev \cite{D} 
coincides with its  resolution  given by
Cattaneo, Felder and Tomassini in  \cite{CFT}. We also give explicitely the identification
of smooth functions, tensorfields and differential operators on $M$ with flat sections
in the corresponding bundles.

Let $M$ be a $d$-dimensional manifold and consider
 the {\bf jet bundle} $E\rightarrow M$ (the bundle of infinite jet of functions)
with fibers $\R [[y^1,\ldots,y^d]]$ (i.e. formal power series in $y\in \R^d$ with
real coefficients) and transition functions induced from the transition functions
of the tangent bundle $TM$. Thus 
\begin{equation}
E=F(M)\times_{\op{Gl}(d,\mathbb{R})}\R [[y^1,\ldots,y^d]]
\end{equation}
where $F(M)$ is the frame bundle. 
 Remark that $E$ can be seen as the formally completed symmetric algebra of the
 cotangent bundle $T^*M$; a section $s\in \Gamma(M,E)$ can be written in the form
$$
s=s(x;y)=\sum_{p=0}^\infty s_{i_1\ldots i_p}(x)y^{i_1}\cdots y^{i_p}
$$Ê
with repeated indices variing from $1$ to $d$, and
where the $s_{i_1\ldots i_p}$ are components of  symmetric covariant tensors on $M.$
This bundle $E$ is denoted  $\mathcal{S}M$ by Dolgushev.\\

  The construction  of a star product on a $d$-dimensional
 Poisson manifold $(M,P)$  given by Cattaneo, Felder
 and Tomassini in \cite{CFT}, using
  a linear torsionfree connection $\nabla$ on the manifold $M$, starts with
the identification of the commutative algebra $C^\infty(M)$ of smooth functions on  $M$
 with the algebra $\mathcal{Z}^0(\Gamma(M,E),D^G)$ of  flat sections  of the jet bundle 
 $E\rightarrow M$ ,
 for the Grothendieck connection $D^G$ (which is constructed using $\nabla$).
 Let us recall this construction.\\

\noindent The exponential map for the connection $\nabla$ gives an identification
\begin{equation}
\exp_x: U\cap T_xM  \rightarrow M ~~y\mapsto \exp_x(y)
\end{equation}
at each point $x$, of the intersection of  the tangent space $T_xM$  with
a  neighborhood $U$ 
of the zero section of the tangent bundle $TM$ 
with a neighborhood of $x$ in $M$. \\Ê\\

\noindent To a function $f\in C^\infty(M)$, one associates 
the section $f_\phi$ of  the jet bundle $E\raa M$ given, for any $x\in M$
by the Taylor expansion 
at $0\in T_xM$ of the pullback $f\circ \exp_x$.
\begin{lemma}
The section $f_\phi$ is given by:
\begin{equation}\label{extaylor}
    f_{\phi}(x;y)=f(x)+\sum_{n>0}\frac{1}{n!} \nabla^n_{i_1\ldots i_n}f(x)\,y^{i_1}\ldots y^{i_n}
    =f(x)+\sum_{n>0}\frac{1}{n!} \nabla^{n,sym}_{i_1\ldots i_n}f(x)\,y^{i_1}\ldots y^{i_n}.
\end{equation}
\end{lemma}
\begin{proof}
In local coordinates $x^i$'s one has 
 $\frac{d}{dt} f(\exp_x ty)=\sum_{k=1}^d(\p_{x^k}f)(\exp_xty)\,\frac{d}{dt} (\exp_xty)^k$ and \\
 $\frac{d^2}{dt^2} f(\exp_xty)=\sum_{k,l}(\p^2_{x^k x^l}f)(\exp_xty)\,\frac{d}{dt} (\exp_xty)^k\,\frac{d}{dt} (\exp_xty)^l+\sum_{k=1}^d(\p_{x^k}f)(\exp_xty)\frac{d^2}{dt^2} (\exp_xty)^k.$
The definition of the exponential map imply that
\begin{equation}\label{eq:geodesics}
 \frac{d^2}{dt^2} (\exp_xty)^k=-\sum_{r,s}\Gamma_{rs}^k(\exp_xty)\,\frac{d}{dt}(\exp_xty)^r\,
 \frac{d}{dt}(\exp_xty)^s
\end{equation}
hence
$\frac{d^2}{dt^2} f(\exp_xty)=\sum_{k,l}(\nabla^2_{kl} f)(\exp_xty)\,\frac{d}{dt} (\exp_xty)^k
\,\frac{d}{dt} (\exp_xty)^l$.
By induction, one gets
$$
\frac{d^n}{dt^n} f(\exp_xty)=\sum_{k_1,\ldots k_n} (\nabla^n_{k_1\ldots k_n} f)(\exp_xty)
\,\frac{d}{dt} (\exp_xty)^{k_1}\ldots\frac{d}{dt} (\exp_xty)^{k_n}
$$
and the result follows at $t=0.$
\end{proof}

\begin{defi}\cite{CF}
The Grothendieck connection $D^G$ on  $E$ is defined by:
\begin{equation}\label{eq:defDG}
D^G_Xs (x;y):=\frac{d}{dt}_{\vert_{t=0}} s(x(t); \exp^{-1}_{x(t)}(\exp_x(y)))
\end{equation}
for  any   curve $t\rightarrow x(t)\in M$
representing $X\in T_xM$ and for any $s\in \Gamma(M,E)$.
It is locally given by
\begin{equation}
D^G_{X}=\sum_{i=1}^{d}X^i \left(\p_{x^i}+\sum_k 
\sum_j\left(\frac{\p \phi_x}{\p y} ^{-1}\right)^k_j \frac{\p \phi^j}{\p x^i}          \p_{y^k}\right)   
\end{equation}
where $\phi_x(y)=\phi(x,y)$ is the Taylor expansion of $\exp_xy$ at $y=0$:
$$
\phi(x,y)^k=x^k+y^k-\frac{1}{2}\sum_{rs}\Gamma_{rs}^k(x)\,y^r\,y^s+
\frac{1}{3!}\sum_{rst}\left(-(\p_{x^r}\Gamma_{st}^k)(x)+2\sum_u\Gamma_{rs}^u(x)\Gamma_{ut}^k(x)\right)\, y^r\, y^s\, y^t
+O(y^4).
$$
\end{defi}
\begin{rem}
$\bullet$ From the definition (\ref{eq:defDG}) it is clear that $D^G$ is flat
$( D^G_X\circ D^G_Y-D^G_Y\circ D^G_X=D^G_{[X,Y]})$.\\
$\bullet$ It is also obvious that $D^G(   f_{\phi})=0 ~  \forall f\in C^\infty(M).$
\end{rem}
\begin{lemma}\cite{CF}
Introducing the operator on $E$-valued forms on $M$
\begin{equation}
\delta=\sum_i dx^i\frac{\p}{\p\,y^i},
\end{equation}
one can write
\begin{equation}
D^G =-\delta +\nabla'+A,
\end{equation}
where
\begin{equation}
\nabla'= \sum_i dx^i\left(\p_{x^i}-\sum_{jk}\Gamma_{ij}^k\,y^j\,\p_{y^k}\right)
\end{equation}
is the covariant derivative  on $E$ associated to $\nabla$ and 
where $A$ is a $1$-form on $M$ with values in the fiberwize vectorfields on $E$, 
\begin{equation}\label{eq:defA}
A(x;y)=:\sum_{ik}dx^i\, A^k_i(x;y)\,\p_{y^k}=\sum_{ik}dx^i  \left( -\frac{1}{3}\sum_{rs}R^k_{ris}(x)y^ry^s+0(y^3)\right) \p_{y^k}.
\end{equation}
\end{lemma}
One extends as usual the operator $D^G$ to the space $\Omega(M,E)$ of $E$-valued forms on $M$:
\begin{equation}
D^G =-\delta +\nabla'+A~~{\rm with}~\nabla'=d-\sum_{ijk}dx^i\Gamma_{ij}^k\,y^j\,\p_{y^k}.
\end{equation}
One introduces the operator $\delta^*=\sum_j y^j\, i(\frac{\p}{\p\,x^j}) 
$ on $\Omega(M,E).$   
Clearly $\left(\delta^*\right)^2=0, \delta^2=0$ and for any   
$\omega\in \Omega^q(M,E_p)$, i.e.  a $q$-form
of degree $p$ in $y$, we have $
(\delta\delta^*+\delta^*\delta)\omega=(p+q)\omega.$\\
Defining, for any  $ \omega \in \Omega^q(M,E_p)$ 
\begin{align*}
\delta^{-1}  \omega  & = \frac{1}{p+q}  \delta^*  \omega~&{\rm when}~p+q\neq 0
\\
                                  &=0   &{\rm when}~p=q= 0
 \end{align*}
we see that any $\delta$-closed $q$-form $\omega$ of degree $p$ in $y$, 
when $p+q>0$, writes uniquely as $\omega=\delta\sigma$ with $\delta^*\sigma=0$; 
$\sigma$ is given by $\sigma=\delta^{-1}\omega$. \\
One proceeds by induction on  the degree in $y$ 
 to see that the cohomoly of $D^G$ is concentrated in degree $0$ and that
any flat section of $E$ is determined by its part of degree $0$ in $y$. 
 Indeed a $q$-form $\omega$
 is $D^G$-closed iff
$\delta\omega=(\nabla'+A)\omega$; this implies
that $\delta\omega_p=0$ for $\omega_p$  the terms of lowest order ($p$) in $y$.
When $p+q>0$ we can write $\omega_p=\delta(\delta^{-1}\omega_p)$ 
and $\omega-D^G(\delta^{-1}\omega_p)$
has terms of lowest order at least $p+1$ in $y$. Remark that given any 
section $s$ of $E$ then $s(x;y=0)$ determines a smooth function $f$ on $M$.
If $D^Gs=0$, then $s-f_\phi$ is still $D^G$ closed. By the above,
Its terms of lowest order  in $y$
must be of the form $\delta\sigma$ hence must vanish since we have a $0$-form.
Hence we have:
\begin{lemma}\cite{CF}
Any section of the jet bundle $s\in\Gamma(E)$ is  the Taylor expansion of the pullback of
a smooth function $f$ on $M$ via the exponential map of the connection $\nabla$
if and only if it is horizontal  for the Grothendieck-connection $D^G$:
\begin{equation}
s=    f_{\phi} ~{\rm for ~a}~ f\in C^\infty(M) \Leftrightarrow s\in\Gamma_{hor}(E):=\{\, s'\in \Gamma(E)
\,\vert\, D^Gs'=0\, \}.
\end{equation}
Furthermore,   the cohomology of $D^G$ is concentrated in degree 0.
 In other word, one obtains a ``Grothendieck-resolution'' of the algebra of smooth functions, i.e. 
 $$H^{\bullet}(\Omega(M,E),D^G)=H^{0}(\Omega(M,E),D^G)=\Gamma_{hor}(E)\cong C^{\infty}(M).$$ 
\end{lemma}
\begin{rem}
When a $q$-form $\omega$ is $D^G$ exact, we have written $\omega=D^G\sigma$
where the tensors defining $\sigma$ are given by universal polynomials
(with no-loop concatenations) in the tensors defining $\omega$, the tensors defining $A$,
the curvature of the connection, and their iterated covariant derivatives.
\end{rem}
\begin{lemma}
The $1$-form $A$ on $M$ with values in the fiberwize vectorfields on $E$
is given by $A(x;y)=:\sum_{ik}dx^i\, A^k_i(x;y)\,\p_{y_k}$
where the $A^k_i$ are universal polynomials given by (no-loop) concatenations of
iterative covariant derivatives of the curvature; they are of the form
\begin{equation}
\sum\left(\nabla_{...}R\right)^{j_1}_{i\cdot\cdot}\left(\nabla_{...}R\right)^{j_2}_{j_1\cdot\cdot} \ldots \left(\nabla_{...}R\right)^{k}_{j_{s-1}\cdot\cdot}y^\cdot\ldots y^\cdot.
\end{equation}
In particular $\delta^{-1}A=0$ since the curvature is skewsymmetric in its first two lower arguments.  The $1$-form $A$ 
is uniquely characterized by the fact that  $\delta^{-1}A=0.$ and the fact that $D^G =-\delta +\nabla'+A$ is flat, i.e;
$\left(D^G\right)^2 =0$ which is equivalent to
\begin{equation}
\delta A= R^{\nabla'}+\nabla' A+\frac{1}{2}[A , A]
\end{equation}
for $\frac{1}{2}[A , A](X,Y):=[A(X),A(Y)]$ and 
$R^{\nabla'}=-\frac{1}{2}R_{ijk}^l dx^i \wedge dx^j y^k\frac{\p}{\p y^l}.$
\end{lemma}
\begin{proof}
Any section $s\in \Gamma(M,E)$ writes
$\sum_{p=0}^\infty s^p_{i_1\ldots i_p}(x)y^{i_1}\cdots y^{i_p}$
with symmetric
$p$-covariant tensors $ s^p_{i_1\ldots i_p}$.
Write  $A(x;y)=\sum_{r \ge 2}dx^i (A^{(r)}(x))^k_{i,j_1\ldots j_r}y^{j_1}\ldots y^{j_r}\p_{y^k}$
with  $(A^{(2)}(x))^k_{i,rs}=-\frac{1}{3}\sum_{rs}R^k_{ris}(x)y^ry^s$.
Then the covariant tensors of  $D^G_Xs$ are given by the symmetrisation of
\begin{equation*}
\left(D^G_Xs\right)^p=-i(X)s^{p+1}+\nabla_Xs^p+\sum_{r=0}^{p-2}(A^{(p-r)}(X))^k\, \p_{y^k}s^{r+1}.
\end{equation*}
The fact that  $D^G(f_\phi)=0~\forall  f \in C^\infty (M)$ implies the expression given in the lemma for $A$.
Indeed, the symmetric tensors defining $f_\phi$ are given by $\frac{1}{p!}\nabla^{p,sym}f$
and we must have
\begin{eqnarray*}
0&=&\left(D^G_{\p_{x_i}}f_\phi\right)^p_{j_1\ldots j_p}
=
-\frac{1}{(p+1)!}(\nabla^{p+1,sym}f)_{ij_1\ldots j_p} +\frac{1}{p!}\left(\nabla (\nabla^{p,sym}f)\right)_{ij_1\ldots j_p}\\
&~& \quad\quad +\sum_{r=0}^{p-2}
\left(A^{(p-r)}(x)\right)^k_{i,j_1\ldots j_{p-r}}\,  \frac{1}{r!}(\nabla^{r+1,sym}f)_{kj_{p-r+1}\ldots j_{p}}
\end{eqnarray*}
with the last terms symmetrized in the $j's$.
The commutation of covariant derivatives of a $q$-form $\omega$ gives
$$
(\nabla^{p+2}f)_{klj_1\ldots j_p}-(\nabla^{p+2}f)_{lkj_1\ldots j_p}=
-\sum_{r=0}^{p}R^s_{klj_r} (\nabla^{p}f )_{j_1\ldots j_{r-1}sj_{r+1}\ldots j_p}
$$
and implies by induction that 
$
(\nabla(\nabla^{p,sym}f))_{ij_1\ldots j_p}-(\nabla^{p+1,sym}f )_{ij_1\ldots j_p}
$
is a universal expression contracting covariant derivatives of the curvature tensor
with lower covariant derivatives of $f$ of the form
$$
\left(\nabla_{ \ldots }R\right)^{t_1}_{i\cdot\cdot}\left(\nabla_{ \ldots }R\right)^{t_2}_{t_1\cdot\cdot} \ldots \left(\nabla_{ \ldots }R\right)^{s}_{t_{s-1}\cdot\cdot} (\nabla^{r+1,sym}f )_{s  \ldots }
$$
with the $j's$ put in a symmetrised way at the $\cdot$'s, and for $0\le r\le p-2$.
Hence the expression for $A$.

Observe that $d \delta+\delta d=0$
and also $\delta \nabla'+\nabla' \delta=0$ since $\nabla$ is torsionfree.
Hence $\left(D^G\right)^2 =0$ iff $-\delta A+ R^{\nabla'}+\nabla' A+\frac{1}{2}[A , A]$
vanishes on all sections of $E$; since it is a $2$-form on $M$ with values in the fiberwize
vectorfields on $E$, this must vanish.
\end{proof}

Dolgushev \cite{D}  gave in a similar spirit a construction for a 
 Kontsevich's formality quasi-isomorphism for a general smooth
 manifold. The construction starts again with a torsionfree linear connexion
 $\nabla$ on $M$. A resolution (called Fedosov's resolution
 in Dolgushev's paper)   of the algebra of functions is given using the complex
 of algebras $(\Omega(M,E),D_F)$
 for a flat connexion (differential) $D_F$  defined by
\begin{equation}
 D_F:=\nabla'-\delta+A
\end{equation}
 where $A$ is a $1$-form on $M$ with values in  the fiberwize
vectorfields on $E$,
 obtained by induction on the order in $y$ by the equation
\begin{equation}
A=\delta^{-1}R^{\nabla'}+\delta^{-1}(\nabla^{'} A+\frac{1}{2}[A,A]).
\end{equation}
This implies that $\delta^{-1}A=0$ and $\delta A=R^{\nabla'}+\nabla^{'} A+\frac{1}{2}[A,A]$
so that $A$ coincides with the $1$-form already considered. Hence
\begin{lemma}
The differential $D^G$ and $D_F$ coincide.
\end{lemma}

 Similarly, Dolgushev defined  a resolution  
 of polydifferential operators and 
 polyvectorfields on $M$ using the complexes 
 $(\Omega(M,\mathcal{D}_{poly}), D_F^{{\mathcal{D}}_{poly}})$ and 
  $(\Omega(M,\mathcal{T}_{poly}), D_F^{{\mathcal{T}}_{poly}})$
  where $\mathcal{T}_{poly}$ is the bundle of formal fiberwize
  polyvectorfields on $E$ and $\mathcal{D}_{poly}$ is the bundle of formal fiberwize
  polydifferential operators on $E$.
 A section of ${\mathcal T}_{poly}^{k}$ is of the form 
\begin{equation}
{\mathcal{F}}(x;y)=\sum_{n=0}^{\infty}{\mathcal{F}}_{i_1\ldots i_n}^{j_1\ldots j_{k+1}}(x)y^{i_1}\ldots y^{i_n}\frac{\p}{\p y^{j_1}}\w\ldots\w\frac{\p}{\p y^{j_{k+1}}},
\end{equation}
 where ${\mathcal{F}}_{i_1\ldots i_n}^{j_1\ldots j_{k+1}}(x)$ are coefficients of  tensors, symmetric  in the covariant indices $i_1,\ldots,i_n$ and antisymmetric in the contravariant indices $j_1,\ldots,j_{k+1}$.
 A section of  ${\mathcal D}_{poly}^{k}$ is of the form 
 \begin{equation}
 {\mathcal{O}}(x;y)=\sum_{n=0}^{\infty}{\mathcal{O}}_{i_1\ldots i_n}^{\alpha_1\ldots \alpha_{k+1}}(x)y^{i_1}\ldots y^{i_n}\frac{\p^{|\alpha_1|}}{\p y^{\alpha_1}}\otimes\ldots\otimes\frac{\p^{|\alpha_{k+1}|}}{\p y^{\alpha_{k+1}}},
 \end{equation} where the $\alpha_l$ are multi-indices and ${\mathcal{O}}_{i_1\ldots i_n}^{\alpha_1\ldots \alpha_{k+1}}(x)$ are coefficients of  tensors symmetric in
 the covariant  indices $i_1,\ldots,i_n.$ and symmetric in each block of $\alpha_i$
 contravariant indices.\\
 
 The spaces $\Omega(M,{\mathcal T}_{poly})$ and $\Omega(M,{\mathcal D}_{poly})$ have a formal fiberwise DGLA structure. Namely, the degree of an element in $\Omega(M,{\mathcal T}_{poly})$ ( resp. $\Omega(M,{\mathcal D}_{poly})$ ) is defined by the sum of the degree of the exterior form and the degree of the polyvector field (resp. the polydifferential operator), the bracket on 
 $\Omega(M,{\mathcal T}_{poly})$ is defined by
$[\omega_1\otimes {\mathcal{F}}_1,\omega_2\otimes {\mathcal{F}}_2]_{SN}:=(-1)^{k_1 q_2} 
\omega_1\w \omega_2 \otimes \left [{\mathcal{F}}_1,{\mathcal{F}}_2\right]_{SN}$ for $\omega_i$
 a $q_i$ form and $ {\mathcal{F}}_i$ a section in ${\mathcal T}_{poly}^{k_i}$
and similarly for  $\Omega(M,{\mathcal D}_{poly})$ using the Gerstenhaber bracket.
The differential on $\Omega(M,{\mathcal T}_{poly})$ is 0 and the differential on $\Omega(M,{\mathcal D}_{poly})$ is defined by $\p:=[m_{pf},.]_G$ where $m_{pf}$ is the fiberwize multiplication of formal power series in $y$ of $E$.

\begin{defi}\cite{D}
The differential $D_F^{{\mathcal{T}}_{poly}}$ is defined on  $\Omega(M,{\mathcal{T}}_{poly})$ by
\begin{equation}
 D_F^{{\mathcal{T}}_{poly}}{\mathcal{F}}:=\nabla^{{\mathcal{T}}_{poly}}{\mathcal{F}}
 -\delta^{{\mathcal{T}}_{poly}} {\mathcal{F}}  +[A,{\mathcal{F}}]_{SN}     
 \end{equation}
 where $\nabla^{{\mathcal{T}}_{poly}}{\mathcal{F}}=d{\mathcal{F}}-\left[\sum_{ijk}dx^i\Gamma_{ij}^k\,y^j\,\p_{y_k},{\mathcal{F}}\right]_{SN}$
 and where
$\delta {\mathcal{F}}=\left[\sum_idx^i\frac{\p}{\p y^i} ,{\mathcal{F}}\right]_{SN}$.
Similarly $D_F^{{\mathcal{D}}_{poly}}$ is defined on  $\Omega(M,{\mathcal{D}}_{poly})$ by
\begin{equation}
  D_F^{{\mathcal{D}}_{poly}}{\mathcal{O}}:=\nabla^{{\mathcal{D}}_{poly}}{\mathcal{O}}
 -\delta^{{\mathcal{D}}_{poly}} {\mathcal{O}}  +[A,{\mathcal{O}}]_{G}     
\end{equation}
with $\nabla^{{\mathcal{D}}_{poly}}$ and $\delta^{{\mathcal{D}}_{poly}}$ defined as above with the
Gerstenhaber bracket.\\
Again the cohomology is concentrated in degree $0$
and a flat section $ {\mathcal{F}}\in {\mathcal{T}}_{poly }$ or $ {\mathcal{O}}\in {\mathcal{D}}_{poly} $
is determined by its terms ${\mathcal{F}}_0$ or ${\mathcal{O}}_0$ of order $0$ in $y$; 
it is defined inductively by
$$
{\mathcal{F}}={\mathcal{F}}_0+\delta^{-1}\left( \nabla^{{\mathcal{T}}_{poly}} {\mathcal{F}}
+[A,{\mathcal{F}}]_{SN}    \right)\quad\quad {\mathcal{O}}={\mathcal{O}}_0+\delta^{-1}\left( \nabla^{{\mathcal{D}}_{poly}} {\mathcal{O}}
+[A,{\mathcal{O}}]_{G} \right).
$$
On the other hand, if $s_1\ldots s_{k+1}$ are sections of $E$,
we have for a ${\mathcal{F}}\in\Gamma(M,{\mathcal{T}}_{poly})$:
\begin{equation}\label{eq:DFtensors}
D_F\left({\mathcal{F}}(s_1,\ldots,s_{k+1})\right)=(D_F^{{\mathcal{T}}_{poly}}{\mathcal{F}})(s_1,\ldots,s_{k+1})
 +{\mathcal{F}}(D_Fs_1,\ldots,s_{k+1})+\cdots+{\mathcal{F}}(s_1,\ldots,D_Fs_{k+1})
\end{equation}
 and similarly for a ${\mathcal{O}}\in\Gamma(M,{\mathcal{D}}_{poly})$.
 \end{defi}
 
\begin{defi}\cite{CFT}
As in  Cattaneo et al. we associate to a polyvector field $F\in T^k_{poly}(M)$
a section $F_{\phi}\in \Gamma(M,{\mathcal{T}}_{poly})$ : for a point $x\in M$
 one considers the Taylor expansion (infinite jet) $F_{\phi}(x;y)$ at $y=0$ of the push-forward $(\exp_{x})_{*}^{-1}F(\exp_xy)$. 
Clearly this definition implies that $X_{\phi}(f_{\phi})=\left(Xf\right)_{\phi}$
so that $F_\phi$ is uniquely determined by the fact that
\begin{equation}\label{eq:def2fphi}
F_\phi(f^1_\phi,\ldots,f^{k+1}_\phi)=\left( F(f^1,\ldots,f^{k+1})\right)_\phi \quad\quad \forall \, f^j\in C^\infty(M).
\end{equation}
Similarly we associate to a differential operator $O\in D^k_{poly}(M)$ a section
$O_\phi \in \Gamma(M,{\mathcal{D}}_{poly})$ determined by the fact that
\begin{equation}
O_\phi(f^1_\phi,\ldots,f^{k+1}_\phi)=\left( O(f^1,\ldots,f^{k+1})\right)_\phi \quad\quad \forall \, f^j\in C^\infty(M).
\end{equation}
\end{defi}
Observe that $D_F^{{\mathcal{T}}_{poly}} F_\phi=0$ by \ref{eq:DFtensors} and \ref{eq:def2fphi} and similarly
$D_F^{{\mathcal{D}}_{poly}} O_\phi=0$, hence we have
\begin{prop}
A section of ${\mathcal T}_{poly}$ is $D_F^{{\mathcal{T}}_{poly}}-$horizontal if and only if is 
a Taylor expansion of a polyvectorfield on $M$, i.e. iff
it is of the form $F_\phi$ for some $F\in T^k_{poly}(M)$.\\
Similarly a section of ${\mathcal D}_{poly}$ is $D_F^{{\mathcal{D}}_{poly}}-$horizontal if and only if is 
of the form $O_\phi$ for some $O\in D^k_{poly}(M)$.\\
The terms in such a flat section are defined by tensors which are universal polynomials (involving concatenations of no-loop type) in the tensors defining the polyvectorfield (or differential operator)
on $M$, the curvature tensor and their iterated covariant derivatives. 
\end{prop}
Observe also that a $D_F$ closed section $s\in \Omega^q(M,E)$
or ${\mathcal{F}}\in \Omega^q(M,{\mathcal{T}}_{poly})$ or
 ${\mathcal{O}}\in \Omega^q(M,{\mathcal{D}}_{poly})$
for $q\ge 1$ is the boundary  of a section defined by tensors which are
given by universal polynomials (involving concatenations of no-loop
type) in the tensors defining the section, the curvature tensor and their iterated
covariant derivatives.

The isomorphisms obtained are isomorphisms of  differential graded Lie algebras.

\section{Construction of a universal star products}\label{section:existence}

The construction of a star product on any Poisson manifold by Cattaneo, Felder and 
Tomassini proceeds as follows: one quantize the identification
of the commutative algebra of smooth functions on $M$ with the algebra of flat sections 
of $E$ in the following way.\\
A deformed algebra stucture on $\Gamma(M,E)[[\nu]]$
 is obtained through fiberwize quantization of the jet bundle using Kontsevich 
 star product on $\R^d$. Precisely, one considers the fiberwize Poisson structure on $E$
 defined by $\P_\phi$ and the fiberwize Kontsevich star product on $\Gamma(M,E)[[\nu]]$:
$$ \sigma*^{\P_\phi}_K\tau=\sigma 
\tau+\sum_{n=1}^{\infty}\frac{\nu^n}{n!}U_n(\P_\phi,\ldots,\P_\phi)(\sigma,\tau).
$$
The operator $D^G_X$ is not a derivation of this deformed product;
 one constructs a  flat connection $D$ which is a derivation of $*^{\P_\phi}_K$.
One defines first
 \begin{equation}
D^1_X=X+\sum_{j=0}^{\infty}\frac{\nu^j}{j!}
U_{j+1}({\hat{X}},\P_\phi,\ldots,\P_\phi)
\end{equation}
where ${\hat{X}}:=D^G_X-X$ is a vertical vectorfield on $E$.
The formality equations imply that $D^1_X$ is a derivation of the star product.
Using the fact that $U_1(\xi)=\xi$ for any vector field $\xi$
and that, for $n\ge 2$,  the maps $U_n(\xi,\alpha_2,\ldots,\alpha_n)=0$ 
if $\xi$ is a linear vector field, we see that
\begin{equation}
D^1_X=D^G_X+\sum_{j=1}^{\infty}\frac{\nu^j}{j!}U_{j+1}
({\hat{\hat{X}}},\P_\phi,\ldots,\P_\phi)
\end{equation}
where ${\hat{\hat{X}}}=\sum_i X^i
(-\p_{y^i}+ \sum_k\, A^k_i(x;y)\,\p_{y^k})$ as defined in equation \ref{eq:defA},
so that it is given by universal polynomials 
(with no-loop type concatenatons) in the tensors defining $X$, the curvature tensor
and their iterated derivatives.
The connection $D^1$ is not flat so one deforms it by
$$
D:=D^1+[\gamma,\cdot]_{*^{\P_\phi}_K}
$$
so that $D$ is flat. The $1$-form $\gamma$ is constructed inductively
using the fact that the cohomology of $D^G$ vanishes.\\
The next point  is to identify  series of  functions on $M$
with the algebra of flat sections of this quantized bundle of algebras
 to define the star product on $M$.\\
 This is done by buildind a map $\rho:\Gamma(M,E)[[\nu]]\rightarrow \Gamma(M,E)[[\nu]]$
 so that $\rho \circ D^G=D \circ \rho$. This map is again constructed by induction
 using the vanishing of the cohomology.\\
 All these points show that the star product constructed in this way is universal.\\

 Similarly, one can show that Dolgushev's formality is universal in the following sense.  Dolgushev constructs a
 $L_\infty$ morphism from the differential graded Lie algebra
of polyvectorfields on $M$ to the differential graded Lie algebra of polydifferential
operators on $M$,
giving  its Taylor coefficients, i.e. a collection of   maps $U^D_j$ associating to $j$ multivectorfields
$F_k $ on $M$ a multidifferential operator $U^D_j(\alpha_1,\ldots,\alpha_j )$. The tensors
defining this operator are given by universal polynomials (involving concatenations
of no-loop type) in the tensors defining the $\alpha_j$'s, the curvature tensor and their iterated
cavariant derivatives.\\
Indeed this $L_\infty$-morphism is obtained in two steps
from  the fiberwize Kontsevich formality from
$\Omega(M,{\mathcal{T}}_{poly})$ to $ \Omega(M,{\mathcal{D}}_{poly})$
building  first a twist which depends only on the curvature and its covariant derivatives,
then building a contraction using the vanishing of the $D_F$ cohomology.

\end{document}